\documentclass[12pt,reqno, letterpaper]{amsart}

\newcommand{\indentalign}{\hspace{0.3in}&\hspace{-0.3in}}
\newcommand{\la}{\langle}
\newcommand{\ra}{\rangle}
\renewcommand{\Re}{\operatorname{Re}}

\newcommand{\sech}{\operatorname{sech}}

\usepackage{amssymb}
\usepackage{amsthm}
\usepackage{amsxtra}
\usepackage{multirow}
\usepackage{graphicx}
\usepackage{ifpdf}
\usepackage{color}

\ifpdf
\DeclareGraphicsExtensions{.pdf}
\else
\DeclareGraphicsExtensions{.pstex,.eps}
\fi

\newtheorem{theorem}{Theorem}
\newtheorem{proposition}[theorem]{Proposition}
\newtheorem{lemma}[theorem]{Lemma}

\numberwithin{equation}{section}
\numberwithin{theorem}{section}

\title{Local ill-posedness of the 1D Zakharov system}
\author{Justin Holmer}
\thanks{The author is partially supported by an NSF postdoctoral fellowship.}
\address{University of California, Berkeley}

\subjclass{primary 35Q55, secondary 35Q51, 35R25}
\keywords{Zakharov system, Cauchy problem, local well-posedness, local ill-posedness}

\begin{document}

\begin{abstract}

Ginibre-Tsutsumi-Velo (1997) proved local well-posedness for the Zakharov system
$$\left\{
\begin{aligned}
&i\partial_tu + \Delta u = nu  \\
&\partial_t^2 n - \Delta n = \Delta |u|^2 \\
&u(x,0)=u_0(x) \\ 
&n(x,0)=n_0(x), \; \partial_tn(x,0)=n_1(x)
\end{aligned}
\right. \qquad
\begin{aligned}
&u=u(x,t)\in \mathbb{C} \\
&n=n(x,t)\in \mathbb{R} \\
&x\in \mathbb{R}^d, \;t\in \mathbb{R}
\end{aligned}
$$
for any dimension $d$, in the inhomogeneous Sobolev spaces $(u,n)\in H^k(\mathbb{R}^d)\times H^s(\mathbb{R}^d)$ for a range of exponents $k$, $s$ depending on $d$.  Here we restrict to dimension $d=1$ and present a few results establishing local ill-posedness for exponent pairs $(k,s)$ outside of the well-posedness regime.  The techniques employed are rooted in the work of Bourgain (1993), Birnir-Kenig-Ponce-Svanstedt-Vega (1996), and Christ-Colliander-Tao (2003) applied to the nonlinear Schr\"odinger equation.
\end{abstract}

\maketitle

\section{Introduction}

\begin{figure}
\label{F:GTV}
\scalebox{0.80}{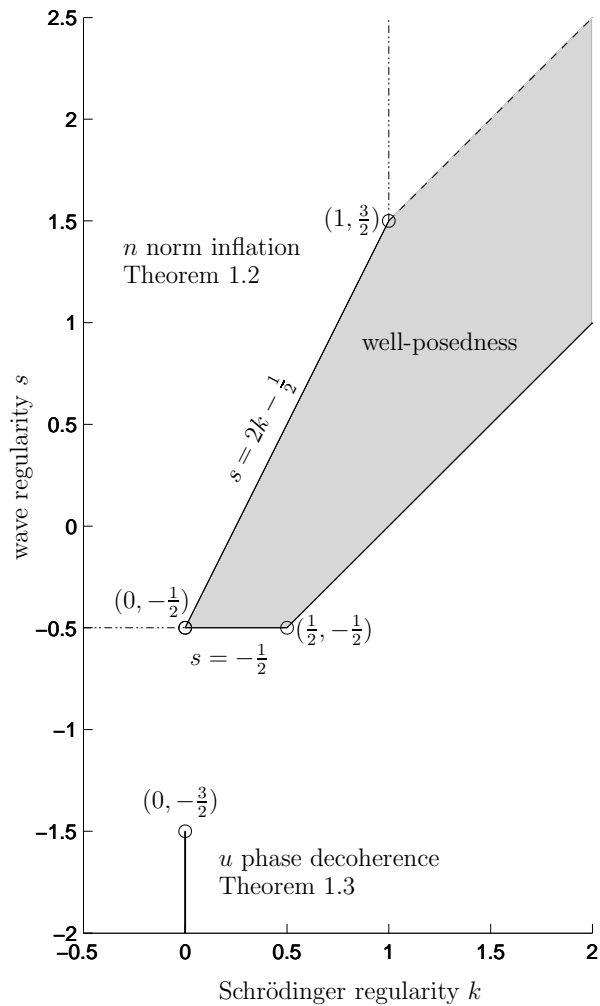}
\caption{The enclosed gray-shaded strip, which extends infinitely to the upper-right, gives the set of pairs $(k,s)$ for which well-posedness has been established by \cite{GTV97} (see Theorem \ref{T:GTVProp1.2}) for $(u_0,n_0,n_1)\in H^k\times H^s\times H^{s-1}$.  Solid lines are included in the well-posedness region, while the dashed line is not.  Theorem \ref{T:northip} provides an ill-posedness result of type ``norm inflation in $n$'' inside the region bounded by the horizontal dotted line $s=-\frac12$, the slanted line $s=2k-\frac12$, and the vertical dotted line $k=1$.  Theorem \ref{T:deepsouthip} provides an ill-posedness result of type ``phase decoherence in $u$'' along the solid vertical line extending down from the point $(0,-\frac32)$.}
\end{figure}

In this paper, we examine the one-dimensional Zakharov system (1D ZS)
$$
\text{1D ZS} \quad \left\{
\begin{aligned}
&i\partial_tu + \partial_x^2 u = nu  \\
&\partial_t^2 n - \partial_x^2 n = \partial_x^2 |u|^2 \\
&u(x,0)=u_0(x) \\ 
&n(x,0)=n_0(x), \; \partial_tn(x,0)=n_1(x)
\end{aligned}
\right. \qquad
\begin{aligned}
&u=u(x,t)\in \mathbb{C} \\
&n=n(x,t)\in \mathbb{R} \\
&x\in \mathbb{R}, \;t\in \mathbb{R}
\end{aligned}
$$
Local well-posedness in the inhomogeneous Sobolev spaces $(u,n)\in H^k(\mathbb{R}) \times H^s(\mathbb{R})$ has been obtained by means of the contraction method in the Bourgain space
$$
\|u\|_{X^S_{k,b_1}} = \left( \iint_{\xi,\tau} \la \xi \ra^{2k} \la \tau + |\xi|^2 \ra^{2b_1} |\hat{u}(\xi,\tau)|^2 \, d\xi \; d\tau \right)^{1/2}
$$
by Bourgain-Colliander \cite{BC96} and Ginibre-Tsutsumi-Velo \cite{GTV97}.\footnote{Actually, these papers consider, more generally, the system in dimensions $d=2,3$ and $d\geq 1$, respectively.}  In the latter paper, the following result is obtained:

\begin{theorem}[\cite{GTV97} Prop.\ 1.2]
\label{T:GTVProp1.2}
\textnormal{1D ZS} is locally well-posed for initial data $(u_0,n_0,n_1) \in H^k\times H^s \times H^{s-1}$ provided 
\begin{align*}
k&\geq 0 & s&\geq -\tfrac{1}{2}\\
-1&\leq s-k < \tfrac{1}{2} &s&\leq 2k-\tfrac{1}{2}
\end{align*}
Specifically:
\begin{enumerate}
\item \emph{Existence}.  $\forall \; R>0$, if $\|u_0\|_{H^k}+ \|n_0\|_{H^s} + \|n_1\|_{H^{s-1}}< R$, then $\exists \; T=T(R)$ and a solution $(u,n)$ to \textnormal{1D ZS} on $[0,T]$ such that
$$
\begin{aligned}
&\|u\|_{C([0,T]; H_x^k)} \leq c\| u_0\|_{H^k} \\
&\|n\|_{C([0,T]; H_x^s)} + \|\partial_t n\|_{C([0,T]; H_x^{s-1})} \leq c\la \|u_0\|_{H^k}\ra ^2 (\|n_0\|_{H^s} + \|n_1\|_{H^{s-1}})
\end{aligned}
$$
and $u\in X_{k,b_1}^S$, where $b_1$ is given by Table \ref{TAB:GTV}.
\item \emph{Uniqueness}.
This solution is unique among solutions $u$ belonging to $C([0,T]; H_x^k)\cap X^S_{k,b_1}$.\footnote{1D ZS can be recast as an integral equation in $u$ alone with $W(n_0,n_1)$ solving \eqref{E:110} appearing as a coefficient.  Then, $n$ can be expressed in terms of $u$ and $W(n_0,n_1)$, and therefore $n$ need not enter into the uniqueness claim.}
\item \emph{Uniform continuity of the data-to-solution map}.  For a fixed $R>0$, taking $T=T(R)$ as above, the map $(u_0,n_0,n_1) \mapsto (u,n, \partial_t n)$ as a map from the $R$-ball in $H^k\times H^s\times H^{s-1}$ to $C([0,T]; H_x^k) \times C([0,T]; H_x^s) \times C([0,T]; H_x^{s-1})$ is uniformly continuous.
\end{enumerate}
\end{theorem}
The region of local well-posedness in this theorem is depicted in Fig.\ \ref{F:GTV}.  We shall outline the \cite{GTV97} proof of Theorem \ref{T:GTVProp1.2} in \S \ref{S:GTV} since the estimates are needed in the proof of Theorem \ref{T:northip} in \S \ref{S:northip}.  

Our goal in this paper is to establish local ill-posedness outside of the well-posedness strip, in particular near the optimal corner $k=0$, $s=-\frac{1}{2}$.  That is, we consider the region (1) $s>2k-\frac{1}{2}$ (above the strip), and (2) $s<-\frac{1}{2}$ (below the strip).  In the first region, the wave data $(n_0,n_1)$ is somewhat smoother than the Schr\"odinger data $u_0$.  As a result, the forcing term $\partial_x^2 |u|^2$ of the wave equation, as time evolves, introduces disturbances that are rougher than the wave data, and the wave solution $n$ does not retain its higher initial regularity.  This is quantified in Theorem \ref{T:northip} below.  In the second region, the Schr\"odinger data $u_0$ is somewhat smoother than the wave data $(n_0,n_1)$.  As a result, the forcing term $nu$ of the Schr\"odinger equation introduces disturbances that are rougher than the Schr\"odinger data, and the Schr\"odinger solution $u$ does not retain its higher initial regularity.  This is quantified in Theorem \ref{T:deepsouthip} and \ref{T:nearsouthip} below.  These simplistic explanations are, at least, accurate for $k>0$.  For $k<0$, there are possibly multiple simultaneous causes for breakdown, although we find that our methods still yield information in this setting.

We will draw upon and suitably modify techniques developed by Birnir-Kenig-Ponce-Svanstedt-Vega \cite{BKPSV96}, Christ-Colliander-Tao \cite{CCT03b}, and Bourgain \cite{Bou93b}, who addressed ill-posedness issues for the nonlinear Schr\"odinger equation.  For a survey of ill-posedness results for nonlinear dispersive equations, see Tzvetkov \cite{T04}.

Our first result demonstrates that the boundary line $s\leq 2k-\frac{1}{2}$ in Theorem \ref{T:GTVProp1.2} is sharp.

\begin{theorem} \quad
\label{T:northip}
Let $0<k<1$ and $s>2k-\frac{1}{2}$ or $k\leq 0$ and $s>-\frac{1}{2}$.  There exists a sequence $\phi_N\in \mathcal{S}$ such that $\|\phi_N\|_{H^k} \leq 1$ for all $N$ and the corresponding solution $(u_N,n_N)$ to \textnormal{1D ZS} on $[0,T]$ with initial data $(\phi_N,0,0)$ satisfies 
\begin{equation}
\label{E:001}
\|n_N(t)\|_{H_x^s} \geq ct N^\alpha \qquad \text{for }0<t\leq T, \; N \geq ct^{-1}
\end{equation}
where $\alpha=\alpha(k,s)>0$.  The time interval $[0,T]$ here is independent of $N$.
\end{theorem}

The form of ill-posedness appearing in Theorem \ref{T:northip} is referred to as ``norm inflation''.  The result is first reduced to the case where $k>0$ and $s$ is just above the line $s=2k-\frac{1}{2}$.  In this case, Theorem \ref{T:GTVProp1.2} applied with $s=2k-\frac{1}{2}$ (the wave initial data is $0$) provides the existence of a solution $(u_N,n_N)$ on a time interval $T$, independent of $N$, with uniform-in-$N$ control on $\|u_N\|_{X_{k,b_1}^S}$.  The estimates of \cite{GTV97} will enable us to show that $u_N$ is comparable to $e^{it\partial_x^2}\phi_N$ in a slightly stronger norm than $X_{k,b_1}^S$ (on this fixed in $N$ time interval) and then Theorem \ref{T:northip} follows from the fact that \eqref{E:001} holds with $n_N=\square^{-1}\partial_x^2|u_N|^2$ replaced by  $\square^{-1}\partial_x^2|e^{it\partial_x^2}\phi_N|^2$, which can be directly verified.\footnote{$w=\square^{-1}f$ is the solution to $\square w =(\partial_t^2-\partial_x^2)w=f$, $w(x,0)=0$, $\partial_t w(x,0)=0$.} The proof is given in \S\ref{S:northip}.

Our second theorem demonstrates lack of uniform continuity of the data-to-solution map, for any $T>0$, as a map from the unit ball in $H^k\times H^s\times H^{s-1}$ to $C([0,T]; H^k)\times C([0,T]; H^s) \times C([0,T]; H^{s-1})$ for $k=0$ and any $s< -\frac{3}{2}$.  We first show that if one issue is ignored, we can, in a manner similar to \cite{BKPSV96}, make use of an explicit soliton class to demonstrate that for any $T>0$ there are two waves, close in amplitude on all of $[0,T]$, initially of the same phase but that slide completely out of phase by time $T$.  This form of ill-posedness is termed ``phase decoherence''.  The soliton class for 1D ZS that we use appears in \cite{Guo89} \cite{Wu94}.  The ``ignored issue'' pertains to low frequencies of $n_0(x)$, and can be resolved by invoking the method of \cite{CCT03b} to construct a ``near soliton'' class offering more flexibility than the exact explicit soliton class in the selection of $n_0(x)$.  This is, however, not straightforward since 1D ZS lacks scaling and Galilean invariance, which was used to manufacture the solution class in \cite{CCT03b}.

\begin{theorem}
\label{T:deepsouthip}
Suppose $k=0$, $s<-\frac{3}{2}$.  Fix any $T>0$ and $\delta>0$.  Then there is a pair of Schwartz class initial data tuples $(u_{0},n_{0},0)$ and $(\tilde{u}_{0}, \tilde{n}_{0}, 0)$ giving rise to solutions $(u,n)$ and $(\tilde{u}, \tilde{n})$ on $[0,T]$ such that the data is of unit size
$$\|u_0\|_{H^k},  \|n_0\|_{H^s} \sim 1, \qquad \|\tilde{u}_0\|_{H^k},  \|\tilde{n}_0\|_{H^s} \sim 1$$
and initially close
$$\| u_0-\tilde{u}_0 \|_{H^k} + \|n_0-\tilde{n}_0\|_{H^s} \leq \delta$$
but the solutions become well-separated by time $T$ in the Schr\"odinger variable
$$\| u(\cdot, t)-\tilde{u}(\cdot, t) \|_{L_{[0,T]}^\infty H_x^k} \sim 1 \,.$$
\end{theorem}

We expect that this result can be extended to all $k\in \mathbb{R}$ and $s<-\frac{3}{2}$, although preliminary efforts were abandoned since the computations became very lengthy and technical.  The proof of Theorem \ref{T:deepsouthip} appears in \S \ref{S:CCT}.

Our final theorem employs a method of Bourgain \cite{Bou93b}. 
\begin{theorem}
\label{T:nearsouthip}
For any $T>0$, the data-to-solution map, as a map from the unit ball in $H^k\times H^{s}\times H^{s-1}$ to $C([0,T]; H^k)\times C([0,T]; H^s) \times C([0,T]; H^{s-1})$ fails to be $C^2$ for $k\in \mathbb{R}$ and $s<-\frac{1}{2}$. 
\end{theorem}
This is a weaker form of ill-posedness than the phase decoherence of Theorem \ref{T:deepsouthip}, although it covers the full region below the well-posedness boundary $s=-\frac{1}{2}$ of \cite{GTV97}. The proof is given in \S\ref{S:Bourgain}.\\

\noindent\textbf{Acknowledgments}.  I would like to thank Jim Colliander for his clear explanation of how to construct counterexamples to bilinear estimates and for other helpful discussion on this topic.  Also, I would like to thank Guixiang Xu for carefully reading \S1--3 of the paper and pointing out numerous misprints and an error.  Finally, I would like to thank the anonymous referee for providing several helpful suggestions for improvement.

\section{The local theory}
\label{S:GTV}

We outline and review the local well-posedness argument in \cite{GTV97} since the estimates will be needed in the proofs of Theorem \ref{T:northip}, \ref{T:nearsouthip}.

Let $[ U(t) u_0 ]\sphat(\xi) = e^{-it\xi^2}\hat{u}_0(\xi)$ and
$$U \ast_R f(\cdot, t) = \int_0^t U(t-t')f(t') \,dt'$$
denote the Schr\"odinger group and Duhamel operators, respectively.  Define the Schr\"odinger Bourgain spaces $X^S_{k,\alpha}$, $Y^S_k$ by the norms
\begin{equation}
\label{E:112}
\begin{aligned}
\|z\|_{X^S_{k,\alpha}} &= \left( \iint_{\xi,\tau} \la \xi \ra^{2k} \la \tau + |\xi|^2 \ra^{2\alpha} |\hat{z}(\xi,\tau)|^2 \, d\xi \; d\tau \right)^{1/2} \\
\|z\|_{Y^S_k} &=  \left( \int_{\xi} \la \xi \ra^{2k} \left(\int_\tau \la \tau + |\xi|^2 \ra^{-1} |\hat{z}(\xi,\tau)| \, d\tau\right)^2 \; d\xi \right)^{1/2} \,.
\end{aligned}
\end{equation}

Consider an initial wave data pair $(n_0, n_1)$.  Split $n_1 = n_{1L}+n_{1H}$ into low and high frequencies\footnote{This decomposition is needed, for otherwise the estimate in Lemma \ref{L:Group}\eqref{I:Groupwave} would have to be modified to have $\|n_1\|_{H^s}$ in place of $\|n_1\|_{H^{s-1}}$ on the right-hand side}, and set $\hat{\nu}(\xi) = \frac{\hat{n}_{1H}(\xi)}{i\xi}$, so that $\partial_x \nu = n_{1H}$.  Let 

$$\begin{aligned}
W_+(n_0,n_1)(x,t) &= \tfrac{1}{2}n_0(x-t) -\tfrac{1}{2}\nu(x-t) + \tfrac{1}{2} \int_{x-t}^x n_{1L}(y)\, dy\\
W_-(n_0,n_1)(x,t) &= \tfrac{1}{2}n_0(x+t) +\tfrac{1}{2}\nu(x+t) + \tfrac{1}{2} \int_x^{x+t} n_{1L}(y)\, dy
\end{aligned}
$$
so that 
\begin{align*}
&(\partial_t\pm\partial_x)W_\pm(n_0,n_1)(x,t) = \tfrac{1}{2}n_{1L}(x) \\
&W_\pm(n_0,n_1)(x,0)=\tfrac{1}{2}n_0(x)\mp \tfrac{1}{2}\nu(x) \,.
\end{align*}
By setting $n = W_+(n_0,n_1)+W_-(n_0,n_1)$, we obtain a solution to the linear homogeneous problem
\begin{equation}
\label{E:110}
\left\{ \begin{aligned}
& \partial_t^2 n - \partial_x^2 n = 0 &\qquad& t,x\in \mathbb{R}\\
& n(x,0) = n_0(x), \; \partial_t n(x,0) = n_1(x) && n=n(t,x)\in \mathbb{R}
\end{aligned}
\right.
\end{equation}
Let
\begin{equation}
\label{E:101}
W_\pm\ast_R f (x,t) = \tfrac{1}{2}\int_0^t f(x\mp s, t-s) \ ds
\end{equation}
so that
$$(\partial_t \pm \partial_x)W_\pm\ast_R f(x,t) = \tfrac{1}{2}f(x,t) \qquad W_\pm f(x,0)= 0 \qquad \partial_tW_\pm f(x,0) = \tfrac{1}{2}f(x,0) \,.$$
It follows that if we set $n=W_-\ast_R f - W_+\ast_R f$, then 
we obtain a solution to the linear inhomogeneous problem
$$
\left\{ \begin{aligned}
& \partial_t^2 n - \partial_x^2 n = \partial_x f &\qquad & t,x\in \mathbb{R}\\
& n(x,0) = 0, \; \partial_t n(x,0) = 0 &&n(x,t)\in \mathbb{R}
\end{aligned}
\right.
$$
Define the one-dimensional reduced wave Bourgain spaces $X^{W\pm}_{s,\alpha}$, $Y^{W\pm}_s$ as
\begin{equation}
\label{E:113}
\begin{aligned}
\|z\|_{X^{W_\pm}_{s,\alpha}} &= \left( \iint_{\xi,\tau} \la \xi \ra^{2s} \la \tau \pm \xi \ra^{2\alpha} |\hat{z}(\xi,\tau)|^2 \, d\xi \; d\tau \right)^{1/2} \\
\|z\|_{Y^{W_\pm}_s} &= \left( \int_\xi \la \xi \ra^{2s} \left( \int_{\tau}  \la \tau \pm \xi \ra^{-1} |\hat{z}(\xi,\tau)| \, d\tau \right)^2 d\xi \right)^{1/2}\,.
\end{aligned}
\end{equation}
Let $\psi(t)=1$ on $[-1,1]$ and $\psi(t)=0$ outside of $[-2,2]$.  Let $\psi_T(t)= \psi(t/T)$, which will serve as a time cutoff for the Bourgain space estimates.  For clarity, we write $\psi_1(t) = \psi(t)$.  We can now recast 1D ZS as
\begin{equation}
\label{E:Zakharovreduced}
\left\{
\begin{aligned}
&i\partial_tu + \partial_x^2 u = (n_++n_-)u  && x\in \mathbb{R}, t\in \mathbb{R}\\
&(\partial_t\pm \partial_x)n_\pm = \mp \tfrac{1}{2}\partial_x |u|^2 +\tfrac{1}{2}n_{1L} 
\end{aligned}
\right.
\end{equation}
where $n=n_++n_-$, which has the integral equation formulation
$$
\begin{aligned}
u(t) &= U(t)u_0 -iU\ast_R[(n_++n_-)u](t) \\
n_\pm(t) &= W_\pm(t)(n_0,n_1) \mp W_\pm\ast_R(\partial_x|u|^2)(t)\,.
\end{aligned}
$$

\begin{lemma}[Group estimates] \label{L:Group} \quad 
\begin{enumerate}
\item \label{I:GroupSch}
\emph{Schr\"odinger}.  $\|\psi_1(t) U(t)u_0\|_{X_{k,b_1}^S} \lesssim \|u_0\|_{H^k}$.
\item \label{I:Groupwave}\emph{1-d Wave}. $\|\psi_1(t) W_\pm(t)(n_0,n_1)\|_{X^{W_\pm}_{s,b}} \lesssim \|n_0\|_{H_x^{s}} + \|n_1\|_{H_x^{s-1}} \,.$
\end{enumerate}
\end{lemma}

\begin{lemma}[Duhamel estimates]
\label{L:Duhamel}
Suppose $T\leq 1$.
\begin{enumerate}
\item \label{I:DuhamelSch}\emph{Schr\"odinger}.  If $0\leq c_1<\frac{1}{2}$, $0\leq b_1$, $b_1+c_1\leq 1$, then $\| \psi_T U\ast_R f \|_{X_{k,b_1}^S} \lesssim T^{1-b_1-c_1}\|f\|_{X_{k,-c_1}^S}$. \\
If $0\leq b_1\leq \frac{1}{2}$, then $\| \psi_T U\ast_R f \|_{X_{k,b_1}^S} \lesssim T^{\frac{1}{2}-b_1}(\|f\|_{X_{k,-\frac{1}{2}}^S\cap Y^S_k})$. \\
$\| U\ast_R f\|_{C(\mathbb{R}_t; H_x^k)} \lesssim \|f\|_{Y^S_k}$.
\item \label{I:Duhamelwave}\emph{1-d Wave}.  If $0\leq c<\frac{1}{2}$, $0\leq b$, $b+c\leq 1$, then $\| \psi_T W_\pm \ast_R f \|_{X_{s,b}^{W\pm}} \lesssim T^{1-b-c}\|f\|_{X_{s,-c}^{W\pm}}$. \\
If $0\leq b\leq \frac{1}{2}$, then $\| \psi_T W_\pm \ast_R f \|_{X_{s,b}^{W\pm}} \lesssim T^{\frac{1}{2}-b}(\|f\|_{X_{s,-\frac{1}{2}}^{W\pm}\cap Y^{W\pm}_s})$.\\
$\|W_\pm \ast_R f\|_{C(\mathbb{R}_t; H_x^s)} \lesssim \|f\|_{Y^{W\pm}_s}$.
\end{enumerate}
\end{lemma}

\begin{lemma}[{\cite{GTV97}} Lemma 4.3/4.5]
\label{L:GTV4.3}
Let $k,s,b,c_1,b_1$ satisfy
\begin{align*}
&s\geq -\tfrac{1}{2} && k\geq 0 && s-k \geq -1 \\
&b,c_1, b_1 > \tfrac{1}{4} && b+c_1 > \tfrac{3}{4} && b+b_1 >\tfrac{3}{4} \\
&s-k \geq -2c_1 
\end{align*}
Then 
$$\|n_\pm u \|_{X^S_{k,-c_1}\cap Y^S_k} \lesssim \|n_\pm \|_{X^{W\pm}_{s,b}} \|u\|_{X^S_{k,b_1}}\,.$$
\end{lemma}

\begin{lemma}[{\cite{GTV97}} Lemma 4.4/4.6]
\label{L:GTV4.4}
Let $k,s,c,b_1$ satisfy
\begin{align*}
& s-2k \leq -\tfrac{1}{2} && k \geq 0 && s-k < \tfrac{1}{2} \\
& c,b_1 > \tfrac{1}{4} && c+b_1 >\tfrac{3}{4} \\ 
& s-k \leq 2b_1-1 && s-k < 2c-\tfrac{1}{2}
\end{align*}
Then 
$$\| \partial_x (u_1 \bar{u}_2) \|_{X^{W\pm}_{s,-c}\cap Y^{W\pm}_s} \lesssim \|u_1\|_{X^S_{k,b_1}}\|u_2\|_{X^S_{k,b_1}} \,.$$
\end{lemma}

To obtain Theorem \ref{T:GTVProp1.2}, fix $0<T<1$, and consider the maps $\Lambda_S$, $\Lambda_{W\pm}$
\begin{align}
\label{E:102} \Lambda_S(u,n_\pm) &= \psi_1Uu_0 + \psi_T U\ast_R[(n_++n_-) u] \\
\label{E:103} \Lambda_{W\pm}(u) &= \psi_1W_\pm(n_0,n_1) \pm \psi_T W_\pm \ast_R(\partial_x |u|^2) \,.
\end{align}
For $T=T(\|u_0\|_{H^k}, \|n_0\|_{H^s}, \|n_1\|_{H^{s-1}})$, a fixed point $(u(t), n_\pm(t)) = (\Lambda_S(u,n_\pm),\Lambda_{W\pm}(u))$ is obtained in $X^S_{k,b_1} \times X^{W\pm}_{s,b}$ satisfying 
\begin{align}
\label{E:104} &\|u\|_{X^S_{k,b_1}} \lesssim \|u_0\|_{H^k} \\
\label{E:105} &\|n\|_{X^{W\pm}_{s,b}} \lesssim \|n_0\|_{H^s} + \|n_1\|_{H^{s-1}} + \|u_0\|_{H^k}^2
\end{align}
by applying Lemmas \ref{L:Group}, \ref{L:Duhamel}, \ref{L:GTV4.3}, \ref{L:GTV4.4} with values for $b_1$, $c_1$, $b$, $c$ given by Table \ref{TAB:GTV}.

\begin{table}
\begin{center}
\renewcommand{\arraystretch}{1.25}
\noindent \begin{tabular}{|c||l|l|}
\hline 
\multirow{2}{*}{$s-k= -1$} & $b_1=\frac{1}{2}-\epsilon$ &  $b = \frac{3}{4}-3\epsilon$  \\
 & $c_1=\frac{1}{2}$  & $c=\frac{1}{4}+2\epsilon$\\
\hline 
\multirow{2}{*}{$-1<s-k< -\frac{1}{2}$} & $b_1=\frac{s-k}{2}+1-\epsilon$  & $b = \frac{3}{4}-2\epsilon$  \\
& $c_1=-\frac{s-k}{2}$ & $c=\frac{1}{4}+\epsilon$ \\
\hline
\multirow{2}{*}{$-\frac{1}{2}\leq s-k \leq 0$} & $b_1=\frac{3}{4}-2\epsilon$ & $b = \frac{3}{4}-2\epsilon$ \\
 & $c_1=\frac{1}{4}+\epsilon$ & $c=\frac{1}{4}+\epsilon$ \\
\hline
\multirow{2}{*}{$0\leq s-k < \frac{1}{2}$} & $b_1=\frac{3}{4}-2\epsilon$ & $b = \frac{3}{4}-\frac{s-k}{2}-2\epsilon$  \\ 
& $c_1=\frac{1}{4}+\epsilon$ & $c=\frac{s-k}{2}+\frac{1}{4}+\epsilon$ \\
\hline
\end{tabular}
\end{center}
\caption{Values of $b_1$, $c_1$, $b$, $c$ meeting the criteria of Lemmas \ref{L:GTV4.3}, \ref{L:GTV4.4} for various intervals of $s-k$.  Note that $b_1+c_1\leq 1-\epsilon$ and $b+c\leq 1-\epsilon$ in order to capture a factor $T^\epsilon$ from Lemma \ref{L:Duhamel}. Also note that $b_1,b>\frac{1}{2}$ and $c_1,c<\frac{1}{2}$ for all cases except $s-k=-1$.}
\label{TAB:GTV}
\end{table}

Consider first the case $s-k>-1$.  We note from Table \ref{TAB:GTV} that $b_1, b>\frac{1}{2}$, and thus we have the Sobolev imbeddings
\begin{equation}
\label{E:107}
\begin{aligned}
&\| u\|_{C(\mathbb{R}_t; H_x^k)} \lesssim \|u\|_{X^S_{k,b_1}} \\
&\|n_\pm\|_{C(\mathbb{R}_t; H_x^s)} \lesssim \|n_\pm \|_{X^{W\pm}_{s,b}}\,.
\end{aligned}
\end{equation}
Also, 
$$\partial_t n(x,t) = \partial_t(n_++n_-)(x,t) = \partial_x(-n_++n_-)(x,t) + n_{1L}(x)$$
and thus
\begin{equation}
\label{E:108}
\| \partial_t n \|_{C(\mathbb{R}_t; H_x^{s-1})} \lesssim \|n_\pm\|_{X^{W\pm}_{s,b}} + \|n_1\|_{H^{s-1}} \,.
\end{equation}
Similar estimates apply to differences of solutions.

Consider now the case $s-k=-1$, where it is necessary to take $b_1<\frac{1}{2}$.  We return to \eqref{E:102} and estimate directly using Lemma \ref{L:Duhamel} to obtain
$$\|u\|_{C(\mathbb{R}_t; H_x^k)} \lesssim \|u_0\|_{H^k} + \|n_\pm u \|_{Y^S_k}$$
and by Lemma \ref{L:GTV4.3}, 
$$\|n_\pm u \|_{Y^S_k} \lesssim \|n_\pm \|_{X^{W\pm}_{k,b}} \|u\|_{X^S_{s,b_1}}$$
where $b_1$, $b$ are as specified in the Table \ref{TAB:GTV}, and  the right-hand side is appropriately bounded by \eqref{E:104}, \eqref{E:105}.  The bounds in \eqref{E:107}, \eqref{E:108} apply in this case since $b>\frac{1}{2}$.
We further note that we can re-estimate $u$ in $X^S_{k,\frac{1}{2}}$ in \eqref{E:102} to obtain
\begin{equation}
\label{E:106}
\|u\|_{X^S_{k,\frac{1}{2}}} \lesssim \|u_0\|_{H^k} + (\|n_0\|_{H^s}+\|n_1\|_{H^{s-1}} + \|u_0\|_{H^k}^2 )\|u_0\|_{H^k} \,.
\end{equation}

\section{Wave norm-inflation for $s> 2k-\frac{1}{2}$}
\label{S:northip}

Here we prove Theorem \ref{T:northip}.  In Steps 1--3, the result will be established for $0<k<\frac{7}{4}$ and $s>2k-\frac{1}{2}$ but with $s$ \textit{near} $2k-\frac{1}{2}$.  In Steps 4--5, the general case of the theorem is reduced to the case considered in Steps 1--3.

\begin{proof}
Let $0<k<1$.  Let 
\begin{align*}
\hat{\phi}_{N,A}(\xi) &= N^{\frac{1}{2}-k}\chi_{[-N-\frac{1}{N},-N]}(\xi) \\
\hat{\phi}_{N,B}(\xi) &= N^{\frac{1}{2}-k}\chi_{[N+1,N+1+\frac{1}{N}]}(\xi)
\end{align*}
Let $\phi_N = \phi_{N,A} + \phi_{N,B}$.
Then $\|\phi_N\|_{H^k} \sim 1$.  A solution to the integral equation
\begin{equation}
\label{E:202}
u_N(t) = \psi_1(t) U(t)\phi_N - i \psi_T(t) U\ast_R \{ [W_+\ast_R(\partial_x|u_N|^2) - W_-\ast_R(\partial_x|u_N|^2)]\cdot u_N\}(t)
\end{equation}
provides a solution to 1D ZS with initial data $(\phi_N, 0,0)$ when $n_N$ is defined in terms of $u_N$ as
\begin{equation}
\label{E:211}
n_N = W_+\ast_R(\partial_x|u_N|^2) - W_-\ast_R(\partial_x|u_N|^2)
\end{equation}
By working with the estimates in Lemmas \ref{L:GTV4.3},  \ref{L:GTV4.4} (taking $s=k-\sigma -\frac{1}{2}$ in the discussion of \S\ref{S:GTV}), we obtain a solution $u_N$ to \eqref{E:202} in $X^S_{k-\sigma,\frac{3}{4}-2\epsilon}$ for $0\leq \sigma \leq k$, on $[0,T]$, where $T=T(\|\phi_N\|_{H^{k-\sigma}})$ (thus independent of $N$) satisfying
\begin{equation}
\label{E:210}
\|u_N\|_{C([0,T]; H_x^{k-\sigma})} \leq \| u_N\|_{X^S_{k-\sigma ,\frac{3}{4}-2\epsilon}} \leq \|\phi_N\|_{H^{k-\sigma}} \sim N^{-\sigma}
\end{equation}

\noindent \textbf{Step 1}.  We show that
\begin{equation}
\label{E:209}
\| [(W_+-W_-)\ast_R \partial_x |U\phi_N|^2](t) \|_{H^s} \sim t N^{s-(2k-\frac{1}{2})} \qquad \text{for }N\gtrsim t^{-1}
\end{equation}
That says that \eqref{E:001} holds provided $u_N(t)$ is replaced by the linear flow $U(t)\phi_N$ in \eqref{E:211}.

To show this, note that in the pairing $U(t)\phi_N \overline{U(t)\phi_N}$, there are 4 combinations $U(t)\phi_{N,j}\overline{U(t)\phi_{N,k}}$, where $j,k\in \{ A,B\}$.  We claim that
\begin{equation}
\label{E:203}
[W_+\ast_R \partial_x (U\phi_{N,A} \overline{U\phi_{N,B}})]\sphat(\xi,t) \; \sim \;i\xi t N^{1-2k}e^{-it\xi} h_1(\xi)
\end{equation}
where $h_1(\xi)$ is the ``triangular step function'' with peak at $\xi=-2N-1-\frac{1}{N}$, of width $\frac{2}{N}$, and of height $\frac{1}{N}$, i.e.
$$h_1(\xi) = \begin{cases} \xi - (-2N-1-\frac{2}{N}) & \text{if }\xi \in [-2N-1-\frac{2}{N},-2N-1-\frac{1}{N}] \\ (-2N-1)-\xi & \text{if }\xi \in [-2N-1-\frac{1}{N}, -2N-1] \end{cases}$$
Here, the symbol $\sim$ means that the difference between the two quantities has $H^s$ norm of lower order in $N$.  
It then follows by taking complex conjugates in \eqref{E:203} that 
\begin{equation}
\label{E:204}
[W_+\ast_R \partial_x (U\phi_{N,B} \overline{U\phi_{N,A}})]\sphat(t,\xi) \sim i\xi t e^{-it\xi} N^{1-2k}h_2(\xi)
\end{equation}
where $h_2(\xi)$ is the ``triangular step function'' centered at $2N+1+\frac{1}{N}$, of width $\frac{2}{N}$, and of height $\frac{1}{N}$, i.e.
$$h_2(\xi) = \begin{cases} \xi-(2N+1) & \text{if }\xi \in [2N+1,2N+1+\frac{1}{N}] \\ (2N+1+\frac{2}{N}) - \xi & \text{if }\xi \in [2N+1+\frac{1}{N}, 2N+1+\frac{2}{N}] \end{cases}$$
Hence 
\begin{equation}
\label{E:208}
\| [W_+\ast_R \partial_x (U\phi_{N,A} \overline{U\phi_{N,B}} + U\phi_{N,B} \overline{U\phi_{N,A}})](t)\|_{H_x^s} \sim t N^{s-(2k-\frac{1}{2})}
\end{equation}
We further claim that the $AA$ and $BB$ interactions for the $W_+$ term are of lower order in $N$, i.e. specifically,
\begin{equation}
\label{E:205}
\| W_+ \ast_R \partial_x(U\phi_{N,j}\overline{U\phi_{N,k}})(t) \|_{H^s} \leq N^{s-(2k-\frac{1}{2})-1} \qquad \text{for }j=k=A \text{ and }j=k=B
\end{equation}
Finally, we claim that all of the interactions $AA$, $AB$, $BA$, and $BB$ for the $W_-$ term are of lower order in $N$, i.e.
\begin{equation}
\label{E:206}
\| [W_- \ast_R \partial_x(U\phi_{N,j}U\phi_{N,k})](t) \|_{H^s} \leq N^{s-(2k-\frac{1}{2})-1} \quad \text{for }j,k \in \{ \, A,B \}
\end{equation}
Combining \eqref{E:208}, \eqref{E:205} \eqref{E:206} establishes \eqref{E:209}.  We begin by proving \eqref{E:203}.  Note that
\begin{align*}
& U(t)\phi_{N,A}(x) = N^{\frac12-k}\int_{\xi_1 \in [-N-\frac{1}{N},-N]} e^{ix\xi_1} e^{-it\xi_1^2} \, d\xi_1 \\
& \overline{U(t)\phi_{N,B}(x)} = N^{\frac12-k}\int_{\xi_2 \in [-N-1-\frac{1}{N}, -N-1]} e^{ix\xi_2} e^{it\xi_2^2} \, d\xi_2
\end{align*}
after the change of variable $\xi_2 \mapsto -\xi_2$ in the second equation.  For the remainder of the computation, $\xi_1$ is restricted to $[-N-\frac{1}{N},-N]$ and $\xi_2$ is restricted to $[-N-1-\frac{1}{N}, -N-1]$. By \eqref{E:101},
\begin{align*}
\indentalign W_+ \ast_R \partial_x (U\phi_{N,A} \overline{U\phi_{N,B}})(t) \\
&= N^{1-2k}\int_{s=0}^t \int_{\xi_1} \int_{\xi_2} i(\xi_1+\xi_2) e^{i(x-s)(\xi_1+\xi_2)} e^{-i(t-s)(\xi_1^2-\xi_2^2)} \, d\xi_1 d\xi_2 \, ds \\
&= N^{1-2k}\int_{\xi_1} \int_{\xi_2} i(\xi_1+\xi_2) e^{ix(\xi_1+\xi_2)} e^{it(\xi_1^2-\xi_2^2)} g(t,\xi_1,\xi_2) \, d\xi_1 d\xi_2
\end{align*}
where
\begin{align*}
g(t,\xi_1,\xi_2) &= \int_{s=0}^t e^{-is(\xi_1+\xi_2)} e^{is(\xi_1^2-\xi_2^2)}  \, ds \\
&= \frac{e^{it(\xi_1+\xi_2)(\xi_1-\xi_2-1)} -1 }{i(\xi_1+\xi_2)(\xi_1-\xi_2-1)}
\end{align*}
Since $\xi_1+\xi_2$ is confined to a $\frac{1}{N}$-sized interval around $-2N-1$ and $\xi_1-\xi_2-1$ is confined to a $\frac{1}{N}$-sized interval around $0$, we have that $(\xi_1+\xi_2)(\xi_1-\xi_2-1)$ is confined to a unit-sized interval around $0$.  By the power series expansion for $e^z$, we have $g(t,\xi_1,\xi_2) \sim t$.
\begin{align*}
\indentalign [W_+ \ast_R \partial_x(U\phi_{N,A} \overline{U\phi_{N,B}})(t)]\sphat(\xi, t) \\
&= N^{1-2k}\int_{\xi_1} \int_{\xi_2} i(\xi_1+\xi_2) \delta(\xi_1+\xi_2-\xi) e^{-it(\xi_1^2-\xi_2^2)} g(t,\xi_1,\xi_2) \, d\xi_1 \, d\xi_2
\end{align*}
Using that $e^{-it(\xi_1^2-\xi_2^2)} = e^{-it(\xi_1-\xi_2-1)(\xi_1+\xi_2)}e^{-it(\xi_1+\xi_2)} \sim e^{-it(\xi_1+\xi_2)}$ and that $g(t,\xi_1,\xi_2) \sim t$, we obtain \eqref{E:203}.  \eqref{E:205} and \eqref{E:206} are proved by a similar computation; we only present the proof of \eqref{E:205} in the case $j=k=A$.  For $t\in [0,T]$, 
\begin{equation}
\label{E:207}
W_+ \ast_R \partial_x(U\phi_{N,A} \overline{U\phi_{N,A}})(t) = \int_{\tau,\xi} \xi \, e^{ix\xi} \frac{e^{it\tau}-e^{-it\xi}}{\tau+\xi} g(\tau,\xi) d\tau d\xi
\end{equation}
where
\begin{align*}
g(\tau,\xi) &= \iint_{\substack {\xi=\xi_1+\xi_2 \\ \tau=\tau_1+\tau_2}} [\psi_1 U\phi_{N,A}]\sphat(\xi_1,\tau_1) [\overline{\psi_1 U{\phi}_{N,A}}]\sphat(\xi_2,\tau_2) \\
&= \iint_{\substack {\xi=\xi_1+\xi_2 \\ \tau=\tau_1+\tau_2}} \hat{\psi}_1(\tau_1+\xi_1^2) \hat{\phi}_{N,A}(\xi_1) \hat{\psi}_1(\tau_2-\xi_2^2) \overline{\hat{\phi}_{N,A}}(-\xi_2)
\end{align*} 
In this integral, $\xi_1$ and $\xi_2$ are each confined to a $\frac{1}{N}$ sized interval around $-N$, forcing $\xi$ to lie in a $\frac{1}{N}$ sized interval around $-2N$.  The $\hat{\psi}_1(\tau_1+\xi_1^2)$ and $\hat{\psi}_1(\tau_2-\xi_2^2)$ factors then (essentially) restrict $\tau_1$ to a unit sized interval around $-N^2$ and restrict $\tau_2$ to a unit sized interval around $N^2$, so that $\tau$ is forced to lie within a unit sized interval around $0$.  Consequently,
$$g(\xi,\tau) \begin{cases} \leq \frac{1}{N} & \text{if }(\xi,\tau) \in [-2N-\frac{2}{N},-2N]\times [-1,1] \\ =0 & \text{otherwise} \end{cases}$$
On the support of $g(\xi,\tau)$, the factor $|\tau+\xi| \sim N$.  From \eqref{E:207},
\begin{align*}
\indentalign \| [W_+ \ast_R \partial_x(U\phi_{N,A} \overline{U\phi_{N,A}})](t)\|_{H^s} \\
&\leq N^{1-2k}\left( \int_\xi |\xi|^2 \la \xi\ra^{2s} \left[ \int_\tau \frac{|g(\tau,\xi)|}{|\tau+\xi|} d\tau \right]^2 d\xi \right)^{1/2} \leq N^{s-(2k-\frac{1}{2})-1}\end{align*}

\noindent \textbf{Step 2}.
Also, on this time interval $[0,T]$ independent of $N$, we claim that 
\begin{equation}
\label{E:201}
\|u_N - \psi_1(t)U(t)\phi_N\|_{X^S_{k+\sigma, b_1}} \leq  \|\phi_N\|_{H^{k'}}^2 \|\phi_N\|_{H^{k+\sigma}} \sim N^{2(k'-k)+\sigma}
\end{equation}
where
\begin{equation}
\label{E:216}
b_1 = \begin{cases} \frac{3}{4} - \frac{k+\sigma}{2} & \text{if }0<k+\sigma<\frac{1}{2} \\ \frac{1}{2} & \text{if }\frac{1}{2}\leq k+\sigma<\frac{5}{2} \end{cases} \qquad 
k' = \begin{cases} 0 & \text{if }0<k+\sigma<\frac{1}{2} \\ \frac{k+\sigma}{2}-\frac{1}{4} & \text{if }\frac{1}{2}\leq k+\sigma<\frac{5}{2} \end{cases}
\end{equation}
Note that $2(k'-k)+\sigma$ will be $<0$ provided $\sigma>0$ is not chosen too large.  This says that $u_N(t)$ is well-approximated by the linear flow $\psi_1(t)U(t)\phi_N$ in the \textit{stronger} norm $X_{k+\sigma}^S$.

We now prove \eqref{E:201}.  From \eqref{E:202},
$$\|u_N - \psi_1 U\phi_N\|_{X^S_{k+\sigma,b_1}} \leq 
\begin{cases} \| (W_\pm \ast_R \partial_x |u_N|^2) \cdot u_N \|_{X^S_{k+\sigma,-c_1}} & \text{if }0\leq k+\sigma \leq \frac{1}{2} \\
\| (W_\pm \ast_R \partial_x |u_N|^2) \cdot u_N \|_{X^S_{k+\sigma,-c_1}\cap Y^S_{k+\sigma}} & \text{if }\frac{1}{2}\leq k+\sigma \leq \frac{5}{2} \end{cases}
$$
for $b_1$ as defined above and 
$$c_1 = \begin{cases} \frac{1}{4} + \frac{k+\sigma}{2} & \text{if }0<k+\sigma<\frac{1}{2} \\ \frac{1}{2} & \text{if }\frac{1}{2}\leq k+\sigma<\frac{5}{2} \end{cases}$$
Following with Lemma \ref{L:GTV4.3},
$$\| u_N-\psi_1U\phi_N\|_{X^S_{k+\sigma,b_1}} \leq \| W_\pm\ast_R \partial_x |u_N|^2 \|_{X^{W\pm}_{s',b}} \|u_N\|_{X^S_{k+\sigma,b_1}}$$
where
$$s' = \begin{cases} -\frac{1}{2} & \text{if }0<k+\sigma<\frac{1}{2} \\ k+\sigma-1 & \text{if }\frac{1}{2}\leq k+\sigma<\frac{5}{2} \end{cases} \qquad b = \begin{cases} \frac{1}{2}-\frac{k+\sigma}{2}+\epsilon & \text{if }0<k+\sigma<\frac{1}{2} \\ \frac{1}{4}+\epsilon & \text{if }\frac{1}{2}\leq k+\sigma<\frac{5}{2} \end{cases}$$
By Lemma \ref{L:GTV4.4},
 $$\| W_\pm \ast_R \partial_x |u_N|^2 \|_{X^{W\pm}_{s',b}} \leq \| \partial_x |u_N|^2\|_{X^{W\pm}_{s',-c}} \leq \|u_N\|_{X^S_{k',b_1'}}^2 $$
where 
\begin{align*}
c&=1-b = \begin{cases} \frac{1}{2}+\frac{k+\sigma}{2}-\epsilon & \text{if }0<k+\sigma<\frac{1}{2} \\ \frac{3}{4}-\epsilon & \text{if }\frac{1}{2}\leq k+\sigma<\frac{5}{2} \end{cases} \\
b_1'&= \begin{cases} \frac{1}{4}+\epsilon & \text{if }0<k+\sigma\leq \frac{1}{2} \\ \frac{k+\sigma}{4}+\frac{1}{8} & \text{if }\frac{1}{2}< k+\sigma<\frac{5}{2} \end{cases}
\end{align*}
and $k'$, $b_1$ are defined above.  Note that $b_1'<\frac{3}{4}-2\epsilon$.  Combining,
\begin{align*}
\indentalign \| u_N - \psi_1 U \phi_N \|_{X^S_{k+\sigma,b_1}} \\
&\leq \|u_N\|_{X^S_{k',b_1'}}^2 \|u_N\|_{X^S_{k+\sigma,b_1}} \\
&\leq \|u_N\|_{X^S_{k',b_1'}}^2\|u_N-\psi_1U\phi_N\|_{X^S_{k+\sigma,b_1}} + \|u_N\|_{X^S_{k',b_1'}}^2\| \psi_1U\phi_N\|_{X^S_{k+\sigma,b_1}}
\end{align*}
By \eqref{E:210},
$$\leq \| \phi_N\|_{H^{k'}}^2\|u_N-\psi_1(t)U(t)\phi_N\|_{X^S_{k,b_1}} + \|\phi_N\|_{H^{k'}}^2\|\phi_N\|_{H^{k+\sigma}}$$
Since $\|\phi_N\|_{H^{k'}} \sim N^{-(k-k')}$, provided $N$ is taken large enough and $k'<k$, \eqref{E:201} will follow.

\noindent \textbf{Step 3}. Here, we establish
$$\|n_N(t)\|_{H^s} \geq tN^{s-(2k-\frac{1}{2})} \quad \text{for }N\geq t^{-1}$$
if $0<k\leq \frac{1}{4}$ and $2k-\frac{1}{2}<s\leq 4k-\frac{1}{2}$, or if $\frac14\leq k<1$ and $2k-\frac12<s\leq \frac43k+\frac16$.

To show this, we note that by \eqref{E:211} and \eqref{E:209}, it suffices to show that 
$$\| W_\pm \ast_R \partial_x(|u_N|^2 - |\psi_1U\phi_N|^2)(t) \|_{H^s} \leq 1$$
Writing
\begin{align*}
\indentalign |u_N|^2 - |\psi_1 U\phi_N|^2 \\
&= |u_N-\psi_1 U\phi_N|^2 + 2\Re \; [(u_N-\psi_1U\phi_N)\overline{\psi_1U\phi_N}]
\end{align*}
we see that it suffices to show that
\begin{align}
& \| [W_\pm \ast_R \partial_x |u_N-\psi_1U\phi_N|^2](t) \|_{H^s} \leq 1 \notag\\
& \| [W_\pm \ast_R \partial_x (u_N-\psi_1U\phi_N) \overline{\psi_1 U\phi_N}](t) \|_{H^s} \leq 1  \label{E:213}\\
& \| [W_\pm \ast_R \partial_x (\psi_1 U\phi_N \cdot \overline{u_N-\psi_1U\phi_N})] (t)\|_{H^s} \leq 1 \notag
\end{align}
We focus on the middle estimate \eqref{E:213}; the other two are handled similarly.  As we describe in detail below, by requiring $s$ to lie sufficiently close to (but above) $2k-\frac{1}{2}$, we can assign $\sigma>0$ such that
\begin{equation}
\label{E:212}
s \begin{cases} \leq 2(k+\sigma) - \frac{1}{2} & \text{if }0<k+\sigma <\frac12 \\ \leq k+\sigma  & \text{if }\frac12\leq k+\sigma<\frac{5}{2} \end{cases}
\end{equation}
and also
\begin{equation}
\label{E:214}
k'+(k+\sigma)\leq 2k
\end{equation}
where $k'$ is given in \eqref{E:216}.  Then proceed to estimate the left-hand side of \eqref{E:213} by Lemma \ref{L:Duhamel}\eqref{I:Duhamelwave} as
$$\|u_N - \psi_1 U \phi_N \|_{X^S_{k+\sigma,b_1}} \|\psi_1 U\phi_N\|_{X^S_{k+\sigma,b_1}}$$
By Step 2 and Lemma \ref{L:Group}\eqref{I:GroupSch},
$$\leq \|\phi_N\|_{H^{k'}}^2 \|\phi_N\|_{H^{k+\sigma}}^2 \sim N^{2(k'-k)}N^{2(k+\sigma-k)}$$
By \eqref{E:214}, it follows that the exponent is $\leq 0$.

We now provide the details assigning $\sigma$ in terms of $k$ and $s$.  The condition \eqref{E:214} is equivalent to the restriction 
\begin{equation}
\label{E:215} \sigma \leq \begin{cases} k & \text{if }k \leq \frac{1}{4} \\ \frac{1}{3}k+\frac{1}{6} & \text{if }\frac{1}{4}\leq k \end{cases}
\end{equation}
The following assignments meet the criteria \eqref{E:215} and \eqref{E:212}.
\begin{itemize}
\item If $0<k\leq \frac{1}{4}$, restrict to $s$ such that $2k-\frac{1}{2}<s\leq 4k-\frac{1}{2}$, and set $\sigma = k$. 
\item If $\frac14\leq k<1$, then restrict to $s$ such that $2k-\frac12<s\leq \frac43k+\frac16$ and set $\sigma=\frac13k+\frac16$.
\end{itemize}

\noindent \textbf{Step 4}. Suppose $0<k<1$ and $s>2k-\frac{1}{2}$.  Let $s'$ be such that $s'\leq s$ and $s'$ meets the restrictions outlined in Step 3 with $s$ replaced by $s'$.  Then by Steps 1--3 (with $s$ replaced by $s'$)
$$\|n_N(t)\|_{H^s} \geq \|n_N(t)\|_{H^{s'}} \geq tN^{s'-(2k-\frac{1}{2})} \quad \text{for }N\geq t^{-1}$$
so we can take $\alpha=s'-(2k-\frac{1}{2})$ in the statement of the theorem.\\

\noindent\textbf{Step 5}.  Next, suppose $k<0$ and $s>-\frac{1}{2}$.  By the reasoning of Step 4, it suffices to restrict to $s<\frac32$.  Set $0<k''< \frac{1}{2}s+\frac{1}{4}$, and note that $s>2k''-\frac{1}{2}$.  Clearly $\|u_N(t)\|_{H^k} \leq \|u_N(t)\|_{H^{k''}}$, so we can just appeal to the conclusion of Steps 1--4 applied with $k$ replaced by $k''$.
\end{proof}

\section{A preliminary analysis for $s\leq -\frac{3}{2}$}

Let $f(x) = \sqrt 2 \sech(x)$, which is the unique positive ground state solution to
\begin{equation}
\label{E:501}
-f + \partial_x^2 f + |f|^2 f =0
\end{equation}
Let $f_\lambda(x) = \lambda f(\lambda x)$ and set
\begin{align*}
u_{\lambda,N}(x,t) &= e^{it(\lambda^2 - N^2)}e^{iNx} \sqrt{1-4N^2} f_\lambda (x-2Nt) \\
n_{\lambda,N}(x,t) &= -|f_\lambda (x-2Nt)|^2
\end{align*}
From \eqref{E:501}, it follows that $(u_{\lambda,N},n_{\lambda,N})$ solves 1D ZS for all $\lambda\in \mathbb{R}$ and $-\frac{1}{2}<N<\frac{1}{2}$.  This is the exact soliton class appearing in \cite{Guo89} and \cite{Wu94}.

Our next goal is to prove Theorem \ref{T:deepsouthip} demonstrating phase decoherence ill-posedness for $k=0$, $s<-\frac{3}{2}$.  We first, however, settle for a partial result (Proposition \ref{P:deepsouthippartial}) using a pair from the above exact explicit soliton class.  We include this result since it is clear and straightforward and exhibits the idea behind the proof of the full result (Theorem \ref{T:deepsouthip}), which is considerably more technical and appears in the next section. 

Define the norm $H^s(|\xi|\geq M)$ as 
$$\| \phi\|_{H^s(|\xi|\geq M)} = \left( \int_{|\xi|\geq M} |\xi|^{2s} |\hat\phi(\xi)|^2 \, d\xi \right)^{1/2}$$
The limitation of the following partial result is the use of $H^s(|\xi| \geq M)$ and $H^{s-1}(|\xi| \geq M)$ norms as opposed to the full $H^s$ and $H^{s-1}$ norms.
\begin{proposition}
\label{P:deepsouthippartial}
Suppse $s\leq -\frac{3}{2}$.  Fix any $T>0$, $\delta>0$.  Then $\exists \; M(\delta)$ sufficiently large and $N(\delta)<\frac{1}{2}$ sufficiently close to $\frac{1}{2}$ so that if 
$$\lambda_1=M, \quad \lambda_2=\sqrt{M^2+\frac{\pi}{2T}}$$
then the solutions are of unit size on $[0,T]$, 
\begin{equation}
\label{E:504}
\begin{gathered}
\|u_{\lambda_j,N}(\cdot, t)\|_{L_x^2} \sim 1 \\
\|n_{\lambda_j,N}(\cdot, t)\|_{H^s(|\xi| \geq M)} \sim 1, \quad  \|\partial_t n_{\lambda_j,N}(\cdot, t)\|_{H^{s-1}(|\xi| \geq M)} \sim 1
\end{gathered}
\end{equation}
and are initially close
\begin{equation}
\label{E:502}
\| u_{\lambda_2,N}(\cdot, 0) -u_{\lambda_1,N}(\cdot, 0) \|_{L^2} \leq \delta
\end{equation}
\begin{equation}
\label{E:503}
\begin{gathered}
\| n_{\lambda_2,N}(\cdot, 0) -n_{\lambda_1,N}(\cdot, 0) \|_{H^s(|\xi|\geq M)} \leq \delta\\
\| \partial_t n_{\lambda_2,N}(\cdot, 0) -\partial_t n_{\lambda_1,N}(\cdot, 0) \|_{H^{s-1}(|\xi|\geq M)} \leq \delta
\end{gathered}
\end{equation}
but become fully separated in the $u$-variable by time $T$,
\begin{equation}
\label{E:505}
\| u_{\lambda_2,N}(\cdot, T) -u_{\lambda_1,N}(\cdot, T) \|_{L^2} \sim 1
\end{equation}
\end{proposition}
\begin{proof}
We will select $M=M(\delta)$ sufficiently large later.  Take $0 \leq N < \frac{1}{2}$ sufficiently close to $\frac{1}{2}$ so that $(1-2N)^{1/2}M^{1/2}=1$.  Then since $N\sim \frac{1}{2}$ we have $\sqrt{1-4N^2} \sim (1-2N)^{1/2}$ and noting that $\lambda_1=M$  and $(1-2N)^{1/2}M^{1/2} =1$ gives
\begin{align*}
\indentalign \|u_{\lambda_2,N}(\cdot,0)- u_{\lambda_1,N}(\cdot, 0)\|_{L^2} = (1-2N)^{1/2} \left\| \hat f\left( \frac{\xi}{\lambda_2} \right) - \hat f\left( \frac{\xi}{\lambda_1} \right) \right\|_{L_\xi^2} \\
&= \left\| \hat f\left( \frac{\lambda_1 \xi}{\lambda_2} \right) - \hat f(\xi) \right\|_{L_\xi^2}
\end{align*}
Take $M$ sufficiently large so that $\lambda_1/\lambda_2$ is sufficiently close to $1$ in order to make the above expression $\leq \delta$.  Thus \eqref{E:502} is established.  Next, we establish \eqref{E:503}.  By the change of variable $\xi \mapsto \lambda_1 \xi$
\begin{align*}
\indentalign \| n_{\lambda_2,N}(\cdot, 0) - n_{\lambda_1,N}(\cdot, 0) \|_{H^s(|\xi|\geq M)}^2  \\
&= \lambda_1^{3+2s} \int_{|\xi| \geq 1} \left| \frac{\lambda_2}{\lambda_1}( f^2)\sphat \left( \frac{\xi\lambda_1}{\lambda_2} \right) - (f^2)\sphat(\xi) \right|^2 |\xi|^{2s} \, d\xi
\end{align*}
Since $s\leq -\frac{3}{2}$ we have $\lambda_1^{3+2s} \leq 1$ and the above difference is made $\leq \delta$ by again taking $M$ sufficiently large.  Also
\begin{align*}
\indentalign \| \partial_t n_{\lambda_2,N}(\cdot, 0) - \partial_t n_{\lambda_1,N}(\cdot, 0) \|_{H^{s-1}(|\xi|\geq M)}^2  \\
&= N^2\lambda_1^{3+2s} \int_{|\xi| \geq 1} \left| \frac{\lambda_2^2}{\lambda_1^2}( f^2\,')\sphat \left( \frac{\xi\lambda_1}{\lambda_2} \right) - (f^2\,')\sphat(\xi) \right|^2 |\xi|^{2(s-1)} \, d\xi
\end{align*}
(Here, the notation $'$ indicates the derivative).  Since $s\leq -\frac{3}{2}$ we have $\lambda_1^{3+2s} \leq 1$ and the above difference is made $\leq \delta$ by again taking $M$ sufficiently large.  The statements \eqref{E:504} are proved by similar change of variable calculations.  The need for the restrictions to $|\xi|\geq M$ in \eqref{E:503} is clear from these calculations.   In fact, one can show that for $s<-\frac12$, we have $\|n_{\lambda, N}(\cdot, 0)\|_{H^s} \sim \lambda$ as $\lambda \to +\infty$ due to the $|\xi|\leq \lambda$ frequency contribution. 

Now we establish \eqref{E:505}.  The key observation here is that while $\lambda_2-\lambda_1$ is very small (as $M\to +\infty$), $\lambda_2^2-\lambda_1^2$ is of fixed size $\pi/(2T)$ and thus $e^{iT(\lambda_2^2-\lambda_1^2)} = i$ is purely imaginary.   Now
$$
\|u_2(\cdot, T) - u_1(\cdot, T)\|_{L^2}^2 = \|u_2(\cdot, T)\|_{L^2}^2 + \|u_1(\cdot, T)\|_{L^2}^2 - 2\Re \int_x u_2(x,T) \overline{u_1(x,T)} \, dx
$$
but the last term on the right-hand side is 
$$ -2\Re e^{iT(\lambda_2^2-\lambda_1^2)}(1-4N^2) \int_x \lambda_2f(\lambda_2x) \lambda_1f(\lambda_1x) \, dx =0$$
which, combined with \eqref{E:504} gives \eqref{E:505}.
\end{proof}

\section{Schr\"odinger phase decoherence for $s<-\frac{3}{2}$}
\label{S:CCT}

Here, we remove the shortcoming of Proposition \ref{P:deepsouthippartial} (high frequency truncated norms $H^s(|\xi|\geq M)$, $H^{s-1}(|\xi|\geq M)$ used instead of $H^s$, $H^{s-1}$) and prove Theorem \ref{T:deepsouthip}.  The soliton class employed in the proof of Proposition \ref{P:deepsouthippartial} involved assigning
$$n(x,t) = -\lambda^2 |f|^2(\lambda(x-2tN))$$
and thus $\hat{n}(\xi,t) = -\lambda (|f|^2)\sphat(\xi/\lambda) e^{-2itN\xi}$.
Replace $|f|^2$ in the definition of $n$ by $g$ defined by $\hat{g}(\xi) = (|f|^2)\sphat(\xi) \chi_{|\xi|\geq 1}(\xi)$ (i.e.\ the restriction to frequencies $\geq 1$) and set
$$\tilde{n}(x,t) = -\lambda^2 g(\lambda (x-2tN))$$
Then 
$$\| \tilde n(\cdot, t)\|_{H^s} + \| \partial_t \tilde n(\cdot, t)\|_{H^{s-1}} \leq 1, \qquad \text{as }\lambda \to +\infty$$
Unfortunately, $(u,\tilde n)$ is no longer a solution to 1D ZS.  We shall thus adapt the method of Christ-Colliander-Tao \cite{CCT03b} to construct a ``near soliton'' class that grants more flexibility in the selection of the wave initial data.  The method procedes by solving a ``small dispersion approximation'' to the equation, and by introducing scaling and phase translation parameters, building the ``near soliton'' class.  The main new obstacle, in comparison to the work of \cite{CCT03b} applied to the nonlinear Schr\"odinger equation, is that 1D ZS does not possess scaling nor the Galilean (phase shift) identity.  We thus need to carry out the small dispersion approximation for a \textit{modified} Zakharov system with the property that when scaling and phase shift operations are performed, the modified Zakharov system is converted into the true Zakharov system.

Consider fixed initial data $(n_0,u_0)$ (to be defined later).  

\textbf{Step 1}. The solution to the small dispersion approximation
$$i\partial_t v = n_0(x) v$$
with $v(x,0)=u_0(x)$ is
$$v(x,t) = e^{-itn_0(x)}u_0(x)$$

\textbf{Step 2}.  For parameters $\lambda \gg 1$, $0<\nu \ll 1$, $-\frac{1}{2}<N<\frac{1}{2}$, consider the initial-value problem for the modified Zakharov system $u=u^{(\lambda,\nu,N)}$, $n_\pm=n_\pm^{(\lambda,\nu,N)}$
\begin{equation}
\label{E:601}
\left\{
\begin{aligned}
&i\partial_t u + \nu^2\partial_x^2 u = 
\begin{aligned}[t]
&\tfrac{1}{2}\left[n_0\left(x+ \frac{\nu(1+2N)}{\lambda}t \right)+ n_0\left(x- \frac{\nu(1-2N)}{\lambda}t\right)\right]u \\
&+ (n_++n_-)u 
\end{aligned}\\
&\frac{\lambda}{(1-2N)\nu} \partial_t n_+ + \partial_x n_+ = -\tfrac{1}{2}(1+2N) \partial_x |u|^2 \\
&\frac{\lambda}{(1+2N)\nu} \partial_t n_- - \partial_x n_- = \tfrac{1}{2}(1-2N) \partial_x |u|^2 \\
&u(x,0)=u_0(x), \qquad n_\pm(x,0)=0
\end{aligned}
\right.
\end{equation}
If $k\geq 1$ and (implicit constants here depend on $\|u_0\|_{H^k}$ and $\|n_0\|_{H^k}$)
\begin{equation}
\label{E:606}
T\lesssim |\ln \nu|, \qquad \lambda \gtrsim \nu^{-5}
\end{equation}
then, on $[0,T]$, we have the two estimates
\begin{align}
\label{E:602} \| u\|_{L_T^\infty H_x^k} &\lesssim \nu^{-1/2} \\
\label{E:603} \|n_\pm \|_{L_T^\infty H_x^{k-1}} &\lesssim \frac{1-4N^2}{\lambda}
\end{align}
To show this, note first that
\begin{align*}
n_+ = -\tfrac{1}{2}(1+2N) \int_0^{t\mu_+} \partial_x |u|^2\left(x-s, t-\frac{s}{\mu_+}\right) ds \\
n_- = \tfrac{1}{2}(1-2N)\int_0^{t\mu_-} \partial_x |u|^2 \left(x+s, t-\frac{s}{\mu_-}\right) \, ds
\end{align*}
where $\mu_+ = \nu(1-2N) /\lambda$ and  $\mu_- = \nu(1+2N)/\lambda$, and thus for $k\geq 1$,
\begin{equation}
\label{E:604}
\|n_\pm \|_{L_T^\infty H_x^{k-1}} \leq  \frac{(1-4N^2)\nu T}{\lambda} \|u\|_{L_T^\infty H_x^k}^2
\end{equation}
By the energy method applied to \eqref{E:601}, we have
\begin{align*}
\indentalign \| \partial_x^ku(T)\|_{L_x^2}^2 - \|\partial_x^k u(0) \|_{L_x^2}^2 \\
&=-\Re i \int_0^T\int_x \partial_x^k[n_0(\cdots)u] \, \overline{\partial_x^ku} \, dxdt -2\Re i \int_0^T\int_x \partial_x^k[n_\pm u ] \, \overline{\partial_x^k u} \, dxdt \\
&= \text{I} + \text{II}
\end{align*}
Term I will be addressed via the Gronwall inequality, while in estimating II we will produce a small coefficient.  
$$|\text{I}| \leq \|n_0\|_{H^k} \int_0^T \|u(t)\|_{H_x^k}^2 \, dt$$
Term II is decomposed as
\begin{align*}
\text{II} &= -\Re i \sum_{\substack{ \alpha+\beta = k \\ \alpha \leq k-1}} c_{\alpha\beta} \int_0^T \int_x \partial_x^\alpha n_\pm \, \partial_x^\beta u \, \overline{\partial_x^k u} \, dxdt - 2\Re i \int_0^T \int_x \partial_x^k n_\pm \, u \, \overline{\partial_x^k u} \,dxdt\\
&= \text{II}_a + \text{II}_b
\end{align*}
From \eqref{E:604},
\begin{align}
\notag |\text{II}_a| &\leq T\|n_\pm \|_{L_T^\infty H_x^{k-1}}\|u\|_{L_T^\infty H_x^k}^2 \\
\label{E:605} &\leq \frac{(1-4N^2)\nu T^2}{\lambda} \|u\|_{L_T^\infty H_x^k}^4
\end{align}
while for $\text{II}_b$, we integrate by parts (here $\approx$ means up to terms bounded similarly to \eqref{E:605})
\begin{align*}
\text{II}_b &= 2\Re i \int_0^T \int_x \partial_x^{k-1}n_\pm \, \partial_x[ u \,\overline{\partial_x^k u} ] \, dxdt \\
&\approx \int_0^T \int_x \partial_x^{k-1}n_\pm \,\partial_x^{k-1}(iu\overline{\partial_x^2 u} - i\bar{u}\partial_x^2u) \, dxdt \\
&= -\frac{1}{\nu^2}  \int_0^T \int_x \partial_x^{k-1}n_\pm \, \partial_x^{k-1}\partial_t |u|^2 \, dxdt \\
&= -\frac{1}{\nu^2}  \int_x \partial_x^{k-1}n_\pm(T) \, \partial_x^{k-1}|u|^2(T) \, dxdt + \frac{1}{\nu^2}  \int_0^T \int_x \partial_t \partial_x^{k-1}n_\pm \, \partial_x^{k-1}|u|^2 \, dxdt \\
&= \text{II}_{b1} + \text{II}_{b2}
\end{align*}
From \eqref{E:604}, we have
$$
| \text{II}_{b1} | \leq \frac{(1-4N^2)T}{\lambda \nu} \|u \|_{L_T^\infty H_x^k}^4$$
From \eqref{E:601}, $\partial_x^{k-1}\partial_t n_\pm = \mp \mu_\pm \partial_x^k n_\pm \mp \frac{1}{2} \mu_\pm(1\pm2N)\partial_x^k|u|^2$, so
\begin{align*}
\text{II}_{b2}&= \pm \frac{\mu_\pm}{\nu^2} \int_0^T\int_x \partial_x^{k-1}n_\pm \, \partial_x^k|u|^2 dxdt\\
\implies |\text{II}_{b2}|&\leq \frac{(1-4N^2)T^2}{\lambda^2} \|u\|_{L_T^\infty H_x^k}^4
\end{align*}
All together, (using $L^2$ conservation as well),
$$\|u(T)\|_{H_x^k}^2 \leq \|u_0\|_{H^k}^2 + \|n_0\|_{H^k}\int_0^T\|u(t)\|_{H_x^k}^2 \, dt+ \epsilon\|u\|_{L_T^\infty H_x^k}^4$$
with 
$$\epsilon = \frac{(1-4N^2)T}{\lambda}\left( \frac{T}{\lambda} + \frac{1}{\nu} + \nu T \right) \leq \frac{1}{\lambda^{1/2}}$$
where the last inequality follows from the assumptions \eqref{E:606}.  By the Gronwall inequality,
\begin{equation}
\label{E:800}
\|u\|_{L_T^\infty H_x^k}^2 \leq e^{\|n_0\|_{H^k}T}(\|u_0\|_{H_x^k}^2 + \lambda^{-1/2}\|u\|_{L_T^\infty H_x^k}^4)
\end{equation}
Provided we have
\begin{equation}
\label{E:801}
\lambda^{-1/2} \lesssim e^{-\|n_0\|_{H^k}T}\|u_0\|_{H^k}^{-2},
\end{equation}
we obtain, from \eqref{E:800} and a continuity in time argument, the bound
\begin{equation}
\label{E:802}
\| u\|_{L_T^\infty H_x^k}^2 \leq 2e^{\|n_0\|_{H^k}T}\|u_0\|_{H^k}^2
\end{equation}
Now, the assumptions \eqref{E:606} imply \eqref{E:800}; and \eqref{E:802} implies \eqref{E:602} by the first of the assumptions in \eqref{E:606}.

\textbf{Step 3}.  With $u=u^{(\lambda,\nu,N)}$ as defined in Step 2, $v$ as defined in Step 1, and \eqref{E:606} satisfied, we claim that 
$$\| u-v \|_{L_T^\infty H_x^k} \lesssim \nu$$
where the implicit constant depends on $\|u_0\|_{H_x^{k+2}}$ and $\|n_0\|_{H_x^{k+2}}$.

For this, we appeal to the result of Step 2 at the level of $k+2$ derivatives, and then apply the energy method to the difference $u-v$ in $H^k$:
\begin{align*}
\indentalign \| \partial_x^k(u-v)(T)\|_{L_x^2}^2 \\
&=
\begin{aligned}[t]
&-2\Re i\nu^2 \int_0^T \int_x \partial_x^{k+2}u \, \overline{\partial_x^k(u-v)} \, dxdt \\
&-2\Re i \int_0^T\int_x \partial_x^k\left[ n_0\left( x \pm \frac{\nu(1\pm 2N)}{\lambda}t\right)u - n_0(x)v\right] \overline{\partial_x^k(u-v)} \, dxdt \\
&-2\Re i\int_0^T \int_x \partial_x^k[n_\pm u] \, \overline{\partial_x^k(u-v)} \, dxdt
\end{aligned}\\
&= \text{I}+\text{II}+\text{III}
\end{align*}
Direct estimates using \eqref{E:602}, \eqref{E:603} give
\begin{align*}
|\text{I}| &\leq T^2\nu^4 \|u\|_{L_T^\infty H_x^{k+2}}^2 + \tfrac{1}{4}\|\partial_x^k(u-v)\|_{L_T^\infty L_x^2}^2 \\
&\leq cT^2\nu^3 + \tfrac{1}{4}\|\partial_x^k(u-v)\|_{L_T^\infty L_x^2}^2
\end{align*}
By  \eqref{E:602} and \eqref{E:603},
\begin{align*}
| \text{III}| &\leq T^2\|n_\pm \|_{L_T^\infty H_x^k}^2 \|u\|_{L_T^\infty H_x^k}^2 + \tfrac{1}{4}\|\partial_x^k(u-v)\|_{L_T^\infty L_x^2}^2\\
&\leq \frac{c(1-4N^2)^2T^2}{\lambda^2\nu} + \tfrac{1}{4}\|\partial_x^k(u-v)\|_{L_T^\infty L_x^2}^2
\end{align*}
By rewriting
\begin{align*}
\indentalign n_0\left(x \pm \frac{\nu(1\pm 2N)}{\lambda}t \right)u - n_0(x)v \\
&= \left[ \int_0^{\pm \frac{\nu (1\pm 2N)t}{\lambda}} \partial_xn_0(x+ s) \, ds \right] u + n_0(x)(u-v)
\end{align*}
term II can be estimated as
$$| \text{II} | \leq 
\begin{aligned}[t]
&\frac{\nu(1\pm 2N)T^2}{\lambda}\|n_0\|_{H^{k+1}}\|u\|_{L_T^\infty H_x^k} \|u-v\|_{L_T^\infty H_x^k} \\
&+ \|n_0\|_{H_x^k} \int_0^T \|(u-v)(t) \|_{H_x^k}^2 \, dt
\end{aligned}
$$
Combining, and applying the Gronwall inequality, we have
$$\|u-v\|_{L_T^\infty H_x^k}^2 \lesssim e^{cT} \Big[ \frac{\nu (1\pm 2N)^2T^4}{\lambda^2} + T^2\nu^3 + \frac{(1-4N^2)^2T^2}{\lambda^2\nu} \Big]$$
The result follows from the assumptions  \eqref{E:606}.

\textbf{Step 4}. For $-\frac{1}{2}<N<\frac{1}{2}$,  set 
\begin{align*}
&u(x,t) = \lambda(1-4N^2)^{1/2} e^{ixN} e^{-itN^2} u^{(\lambda,\nu,N)}(\lambda\nu (x-2tN), \lambda^2 t) \\
&n_\pm(x,t) = \lambda^2 n_\pm^{(\lambda,\nu,N)}(\lambda \nu (x-2tN), \lambda^2 t)
\end{align*}
Then $(u,n)$ solves 1D ZS, with 
$$n(x,t)=n_+(x,t)+n_-(x,t)+\tfrac{1}{2}\lambda^2n_0(\lambda\nu(x+t))+\tfrac{1}{2}\lambda^2n_0(\lambda\nu(x-t))$$ 
and initial data $u(x,0)=\lambda u_0(\lambda \nu x)$, $n(x,0) = \lambda^2 n_0(\lambda \nu x)$, $\partial_t n(x,0)=0$.

Consider $0< \nu \ll 1$ and $\lambda \gg 1$.  Since $\|u^{(\lambda, \nu , N)}(x,t)\|_{L_x^2} = \|u_0\|_{L_x^2}$ for all $t$, we have by change of variable
$$\|u(x,t)\|_{L_x^2} = \frac{\lambda^{1/2}(1-4N^2)^{1/2}}{\nu^{1/2}} \|u_0\|_{L_x^2}$$
Also, if $\hat n_0(\xi) = 0$ for $|\xi|\leq 1$ and $\lambda \nu \geq 1$, then another change of variable gives
$$\|n(x,0)\|_{H_x^s} \leq \lambda^{\frac{3}{2}+s}\nu^{s-\frac{1}{2}}\|n_0\|_{H_x^s}$$
If $s<-\frac{3}{2}$ and $\lambda$ and $\nu$ satisfy
\begin{equation}
\label{E:611}
\lambda  \geq \nu^{-\alpha}, \qquad \text{with }\alpha = \max\left( \frac{ \frac{1}{2}-s}{-s-\frac{3}{2}},\, 5\right)
\end{equation}
then $\|n(x,0)\|_{H_x^s} \leq \|n_0\|_{H^s}$. 

\textbf{Step 5}.  Fix $M \gg 1$ and $0\leq \nu \ll 1$, to be chosen momentarily.  In terms of $M$ and $\nu$, define the following quantities:  Let $T=|\ln \nu|/M^2$,  and set
$$\lambda_1=M, \quad \lambda_2 = \sqrt{\frac{\pi}{2T}+M^2}$$
We note that 
\begin{equation}
\label{E:612}
e^{iT(\lambda_2^2-\lambda_1^2)} = i
\end{equation}
is purely imaginary.  Note further that
\begin{equation}
\label{E:608}
\frac{\lambda_2}{\lambda_1} = \sqrt{ \frac{\pi}{2|\ln \nu|} + 1} \to 1, \quad \text{as }\nu \to 0, \quad \text{independently of }M
\end{equation}
Take $N$ sufficiently close to $\frac{1}{2}$ so that 
\begin{equation}
\label{E:610}
\frac{(1-2N)M}{\nu}=1
\end{equation}
Take $u_0(x)\in \mathcal{S}(\mathbb{R})$ such that $u_0(x)=1$ for $-1 \leq x \leq 1$, and $n_0(x)$ to be the smooth function given on the Fourier side as
$$\hat{n}_0(\xi) = 
\left\{
\begin{aligned}
&0 && \text{if }|\xi|\leq 2\\
&\frac{\pi}{2} && \text{if }2\leq |\xi| \leq 4\\
&0 && \text{if }|\xi| \geq 4
\end{aligned}
\right.
$$
so that in fact 
$$n_0(x) = \frac{\cos 3x \sin x}{x}$$
Now consider the solutions $u^{(\lambda_1,\nu,N)}$ and $u^{(\lambda_2,\nu,N)}$ of the modified Zakharov system given at the beginning of Step 2, both in terms of the data $u_0(x)$ and $n_0(x)$.  Define, as in Step 4, the 1D ZS solution $(u_1,n_1)$ in terms of $u^{(\lambda_1,\nu,N)}$ and $(u_2,n_2)$ in terms of $u^{(\lambda_2,\nu,N)}$.  By the comments at the end of Step 4, 
\begin{gather*}
\|u_j(x,t)\|_{L_x^2} =1\\
\|n_j(x,0)\|_{H_x^s} + \|\partial_t n_j(x,0) \|_{H_x^{s-1}} \leq 1
\end{gather*}
where, in order to meet condition \eqref{E:611}, we need $M=M(\nu) \geq \nu^{-\alpha}$.
By a change of variable and \eqref{E:610},
$$
\| u_2(x,t) - u_1(x,t) \|_{L_x^2} =  \left\| \frac{\lambda_2}{\lambda_1} u^{(\lambda_2,\nu,N)}\left( \frac{\lambda_2 x}{\lambda_1}, \lambda_2^2t \right) - u^{(\lambda_1,\nu,N)}(x,\lambda_1^2 t) \right\|_{L_x^2}
$$
By \eqref{E:608} and the fact that $\|u^{(\lambda,\nu,N)}(x,t)\|_{L_x^2} = \|u_0\|_{L_x^2}$ for all $t$ and uniformly in all the parameters, we can take $\nu=\nu(\delta)>0$ sufficiently small so that 
$$\| u_2(x,t) - u_1(x,t) \|_{L_x^2} = \| u^{(\lambda_2,\nu,N)}( x, \lambda_2^2t ) - u^{(\lambda_1,\nu,N)}(x,\lambda_1^2 t) \|_{L_x^2} + \mathcal{O}(\delta)
$$
By the results of Step 1 and 3, again taking $\nu=\nu(\delta)$ sufficiently small, if $0 \leq t \leq |\ln \nu|/M^2$, then
\begin{equation}
\label{E:609}
\| u_2(x,t) - u_1(x,t) \|_{L_x^2} = \| (e^{i(\lambda_2^2-\lambda_1^2)t \, n_0(x)} - 1)u_0(x) \|_{L_x^2} + \mathcal{O}(\delta)
\end{equation}
Here, the first condition of \eqref{E:606} is met, since the $T$ appearing there is our $\lambda_1^2t \sim \lambda_2^2t \lesssim |\ln v|$.  We see trivally from \eqref{E:609} that
$$\| u_2(x,0) - u_1(x,0) \|_{L_x^2} \lesssim \delta$$
But by \eqref{E:612} and \eqref{E:609} and the choice of $u_0(x)$ and $n_0(x)$,
$$\| u_2(x,T) - u_1(x,T) \|_{L_{[0,T]}^\infty L_x^2} = \mathcal{O}(1)$$
We further note that
$$T= \frac{|\ln \nu|}{M^2} \leq |\ln \nu| \nu^4 \to 0 \qquad \text{as }\nu \to 0$$
and therefore we can accomodate an arbitrarily small preselected time, $T$ as in the statement of the theorem.

\section{The Schr\"odinger flow map is not $C^2$ for $s<-\frac{1}{2}$}
\label{S:Bourgain}

In this section, we give the proof of Theorem \ref{T:nearsouthip}.  For fixed $H^\infty$ data $(u_0,n_0,n_1)$, to be specified later, and a parameter $\gamma \in \mathbb{R}$, consider initial data $(u\big|_{t=0},n\big|_{t=0},\partial_tn\big|_{t=0})=(\gamma u_0, \gamma n_0, \gamma n_1)$ and corresponding 1D ZS solutions $(u,n)=(u_\gamma,n_\gamma)$.  Clearly
\begin{equation}
\label{E:701}
u\big|_{\gamma=0}=0, \quad \partial_xu\big|_{\gamma=0}=0, \quad \partial_x^2u\big|_{\gamma=0}=0, \quad n\big|_{\gamma=0}=0, \quad \partial_xn\big|_{\gamma=0}=0
\end{equation}
The solution, written in integral equation form, is:
\begin{align*}
&u(t) = U(t)(\gamma u_0) -i \int_0^t U(t-t')[(un)(t')]\, dt' \\
&n(t) = W(t)(\gamma n_0,\gamma n_1) \pm \tfrac{1}{2} \int_0^t \partial_x |u|^2(x\pm s,t-s) \, ds
\end{align*}
from which it follows that
\begin{equation}
\label{E:702}
\begin{aligned}
&\partial_\gamma u(t) = U(t)u_0 -i \int_0^t U(t-t')[(\partial_\gamma u \, n + u \,\partial_\gamma n)(t')]\, dt' \\
&\partial_\gamma n(t) = W(t)(n_0,n_1) \pm \tfrac{1}{2}\int_0^t \partial_x (\partial_\gamma u \, \bar{u} + u \overline{\partial_\gamma u})(x\pm s,t-s) \, ds
\end{aligned}
\end{equation}
By \eqref{E:701},
\begin{equation}
\label{E:703}
\partial_\gamma u\big|_{\gamma=0} = Uu_0, \quad \partial_\gamma n\big|_{\gamma=0}=W(n_0,n_1)
\end{equation}
By applying $\partial_x$ to \eqref{E:702} and again appealing to \eqref{E:701}, we get
\begin{equation}
\label{E:704}
\partial_x \partial_\gamma u\big|_{\gamma=0}=\partial_xUu_0, \quad \partial_\gamma\partial_x n\big|_{\gamma=0}=\partial_xW(n_0,n_1)
\end{equation}
By applying $\partial_\gamma$ to \eqref{E:702}, we obtain
\begin{align*}
&\partial_\gamma^2 u(t) = -i \int_0^t U(t-t')[(\partial_\gamma^2 u \, n +2\partial_\gamma u \, \partial_\gamma n + u \,\partial_\gamma^2 n)(t')]\, dt' \\
&\partial_\gamma^2 n(t) =  \pm \int_0^t \partial_x (\partial_\gamma^2 u \, \bar{u} +2|\partial_\gamma u|^2+ u \overline{\partial_\gamma^2 u})(x\pm s,t-s) \, ds
\end{align*}
from which we find, together with \eqref{E:703}\eqref{E:704}, that
\begin{align*}
&\partial_\gamma^2 u \big|_{\gamma=0}(t) = -2i \int_0^t U(t-t')[ Uu_0(t') \; W(n_0,n_1)(t')]\, dt'\\
&\partial_\gamma^2 n \big|_{\gamma=0}(t) = \pm \int_0^t \partial_x|U(t')u_0|^2(x\mp s, t-s) \, ds
\end{align*}

Let $X=H^k\times H^{s}\times H^{s-1}$, $Y=H^k\times H^s$.  Fix $t>0$, and let $F:X\to Y$ be the solution map $F(u_0,n_0,n_1) = (u(t),n(t))$.
Let $G:\mathbb{R}\to X$ be given by $G(\gamma)=(\gamma u_0,\gamma n_0, \gamma n_1)$.  Let $H(\gamma)=F\circ G(\gamma)$ so that $H:\mathbb{R} \to Y$.  Then (here $\mathcal{L}(A;B)$ denotes a linear map $A\to B$)
$$\underbrace{DH(\gamma)}_{\mathcal{L}(\mathbb{R};Y)} = \underbrace{DF(G(\gamma))}_{\mathcal{L}(X;Y)}\circ \underbrace{DG(\gamma)}_{\mathcal{L}(\mathbb{R};X)}$$
Also
$$
D^2H(\gamma)=\underbrace{D^2F(G(\gamma))}_{\mathcal{B}(X\times X;Y)}\circ(\underbrace{DG(\gamma)}_{\mathcal{L}(\mathbb{R};X)},\underbrace{DG(\gamma)}_{\mathcal{L}(\mathbb{R};X)})+\underbrace{DF(G(\gamma))}_{\mathcal{L}(X;Y)}\circ \underbrace{D^2G(\gamma)}_{\mathcal{B}(\mathbb{R}\times \mathbb{R}; X)}
$$
and since $D^2G(\gamma)=0$
$$D^2H(\gamma)(\alpha_1,\alpha_2)=D^2F(G(\gamma))((\alpha_1 u_0, \alpha_1 n_0, \alpha_1 n_1), \;(\alpha_2 u_0, \alpha_2 n_0, \alpha_2 n_1)) $$
Hence
$$D^2F(0)((u_0,n_0,n_1),(u_0,n_0,n_1)) = D^2H(0)(1,1)= (\partial_\gamma^2 u\big|_{\gamma=0}(t), \partial_\gamma^2 n\big|_{\gamma=0}(t))$$
We now note how to prescribe an appropriate sequence $(u_{N,0},n_{N,0},n_{N,1})$ (indexed by $N$) to show that $D^2F(0)\in \mathcal{B}(X\times X; Y)$ is not a bounded (continuous) bilinear map in the two cases (1) $s<-\frac{1}{2}$ and (2) $s>2k-\frac{1}{2}$.  
\begin{itemize}
\item If $s<-\frac{1}{2}$, 
$$\hat{u}_{N,0}(\xi) = N^{\frac{1}{2}-k}\chi_{[-N-\frac{1}{N},-N]}(\xi)$$
$$\hat{n}_{N,0}(\xi) = N^{\frac{1}{2}-s}\chi_{[2N-1,2N-1+\frac{1}{N}]}(\xi) + N^{\frac{1}{2}-s}\chi_{[-2N+1-\frac{1}{N},-2N+1]}(\xi)$$
and $n_{N,1}=0$, then $(u_{N,0},n_{N,0},n_{N,1})$ is a sequence  such that $\|(u_{N,0},n_{N,0},n_{N,1})\|_X \sim 1$ but
$$\|\partial_\gamma^2u\big|_{\gamma=0}(t)\|_{H^k} \geq c(t) N^{-s-\frac{1}{2}}$$ We note that the second term in the definition of $\hat n_{N,0}(\xi)$ is included solely to make $n_0(x)$ real.
\item If $s>2k-\frac{1}{2}$ and we set
$$\hat{u}_{N,0}(\xi) = N^{\frac{1}{2}-k}(\chi_{[-N-\frac{1}{N},-N]}(\xi) + \chi_{[N+1,N+1+\frac{1}{N}]}(\xi))$$
and $n_{N,0}=0$, $n_{N,1}=0$,  then $(u_{N,0},n_{N,0},n_{N,1})$ is a sequence  such that $\|(u_{N,0},n_{N,0},n_{N,1})\|_X \sim 1$ but
$$\|\partial_\gamma^2n\big|_{\gamma=0}(t)\|_{H^s} \geq c(t) N^{s-(2k-\frac{1}{2})}$$
\end{itemize}
The second was considered in \S \ref{S:northip} as part of the proof of Theorem \ref{T:northip}, and thus reproduces a weaker version of that result.  We now carry out a proof of the first case to establish Theorem \ref{T:nearsouthip}.
Since
\begin{equation}
\label{E:705}
[\partial_\gamma^2u\Big|_{\gamma=0}(t)]\sphat(\xi) = \int_0^t e^{-i(t-t')\xi^2}[Uu_{N,0}(t')W(n_{N,0},0)(t')]\sphat(\xi) \, dt'
\end{equation}
we need to examine
\begin{align*}
\indentalign [Uu_{N,0}(t')W(n_{N,0},0)(t')]\sphat(\xi) = \int_{\xi_1} e^{-it'\xi_1^2}\hat{u}_{N,0}(\xi_1) \cos(t'(\xi-\xi_1)) \hat n_{N,0}(\xi-\xi_1) \, d\xi_1 \\
&= (e^{-it'N^2}\cos(t'(2N-1))+ \mathcal{O}(t')) \int_{\xi_1} \hat{u}_{N,0}(\xi_1) \hat n_{N,0}(\xi-\xi_1) \, d\xi_1 
\end{align*}
by the support properties of $u_{N,0}$ and $n_{N,0}$.  Directly evaluating the convolution gives
$$= (e^{-it'(N-1)^2}+e^{-it'(N^2+2N-1)}+\mathcal{O}(t')) N^{1-k-s}g(\xi)$$
where $g(\xi)=g_1(\xi)+g_2(\xi)$ consists of two triangular step functions, each of height $1/N$ and width $2/N$, centered at $N-1$ and $-3N+1$, respectively.  Specifically,
$$g_1(\xi) = 
\left\{
\begin{aligned}
&\tfrac{1}{N} - |\xi-(N-1)| && \text{if }|\xi-(N-1)| \leq \tfrac{1}{N} \\
&0 && \text{otherwise}
\end{aligned}
\right.
$$
$$
g_2(\xi) = 
\left\{
\begin{aligned}
&\tfrac{1}{N} - |\xi-(-3N+1)| && \text{if }|\xi-(-3N+1)| \leq \tfrac{1}{N}\\
&0 && \text{otherwise}
\end{aligned}
\right.
$$
We have by the support properties of $g_1(\xi)$ and $g_2(\xi)$, and \eqref{E:705}
\begin{align*}
\indentalign [\partial_\gamma^2 u\Big|_{\gamma=0}(t)]\sphat(\xi)\\
&= 
\begin{aligned}[t]
&+N^{1-k-s} g_1(\xi) \int_0^t e^{-i(t-t')(N-1)^2}(e^{-it'(N-1)^2} + e^{-it'(N^2+2N-1)}) \, dt' \\
&+N^{1-k-s} g_2(\xi)\int_0^t e^{-i(t-t')(-3N+1)^2}(e^{-it'(N-1)^2} + e^{-it'(N^2+2N-1)}) \, dt'\\
&+N^{1-k-s}g(\xi)\mathcal{O}(t^2)
\end{aligned}
\end{align*}
Evaluating each component separately gives
$$[\partial_\gamma^2 u\Big|_{\gamma=0}(t)]\sphat(\xi) = N^{1-k-s}g_1(\xi)e^{-it(N-1)^2}t +  N^{1-k-s}g(\xi)(\mathcal{O}(t^2) + \mathcal{O}(N^{-1}))$$
Thus, provided $t$ is chosen small and $N$ sufficiently large, the first term is pointwise dominant, giving
$$\|\partial_\gamma^2 u\big|_{\gamma=0}(t)\|_{H^k} \geq tN^{-\frac{1}{2}-s}$$
completing the proof.

\newcommand{\etalchar}[1]{$^{#1}$}
\providecommand{\bysame}{\leavevmode\hbox to3em{\hrulefill}\thinspace}
\providecommand{\MR}{\relax\ifhmode\unskip\space\fi MR }
\providecommand{\MRhref}[2]{%
  \href{http://www.ams.org/mathscinet-getitem?mr=#1}{#2}
}
\providecommand{\href}[2]{#2}

\end{document}